\documentclass[12pt]{amsart}
\usepackage{amssymb}
\usepackage{amsmath}
\usepackage{mathabx}
\usepackage{amsthm}
\usepackage{amstext}
\usepackage{amsopn}
\usepackage{mathrsfs}
\usepackage{latexsym}
\usepackage{bm}
\DeclareMathAlphabet{\mathpzc}{OT1}{pzc}{m}{it}
\allowdisplaybreaks
\newtheorem{Definition}{Definition}[section]
\newtheorem{Proposition}{Proposition}[section]
\newtheorem{Lemma}{Lemma}[section]
\newtheorem{Theorem}{Theorem}[section]
\newtheorem{Corollary}{Corollary}[section]
\newtheorem{Remark}{Remark}[section]
\newtheorem{Example}{Example}[section]

\begin{document}
\bibliographystyle{plain}
\footnotetext{
\emph{2010 Mathematics Subject Classification}: 46L53, 05C76\\
\emph{Key words and phrases:}
monotone independence, conditionally monotone independence, tensor independence, comb product of graphs, c-comb product of graphs}
\title[Conditionally monotone independence]
{Conditionally monotone independence and the associated products of graphs}
\author[R. Lenczewski]{Romuald Lenczewski}
\address{Romuald Lenczewski \newline
Katedra Matematyki, Politechnika Wroc\l{}awska, \newline
Wybrze\.{z}e Wyspia\'{n}skiego 27, 50-370 Wroc{\l}aw, Poland}
\email{Romuald.Lenczewski@pwr.edu.pl}
\begin{abstract}
We reduce the conditionally monotone (c-monotone) independence of Hasebe to 
tensor independence on suitably constructed larger algebras. For that purpose, we use the approach developed for a reduction of 
similar type for boolean, free and monotone independences.
We apply the tensor product realization of c-monotone random variables to introduce the {\it c-comb} ({\it loop}) {\it 
product of birooted graphs}, 
a generalization of the comb (loop) product of rooted graphs, and we show that it is related to the c-monotone additive (multiplicative) convolution of distributions.
\end{abstract}
\maketitle

\section{Introduction}
There are several notions of noncommutative independence, including freeness of Voiculescu [20], related to the free product of groups, boolean independence that can be traced back to the so-called regular free product of groups studied by Bo\.{z}ejko [5] and monotone independence of Muraki [17] and Lu [16]. In noncommutative probability theory, classical independence is known under the name of tensor independence since it is related to the tensor product of algebras.

The axiomatic theory distinguishes the notions mentioned above as those which satisfy certain natural axioms. The early version of the axiomatic approach [3] said that only free, tensor and boolean independences were appropriate notions. Inspired by this work, we demonstrated that all these notions were included in the framework of tensor independence [11,12]. More precisely, 
we showed that boolean independence as well as free independence can be `reduced' to tensor independence by going to tensor products of extended algebras (*-algebras) and states on these algebras. Then, Muraki [17] and Lu [16] introduced the monotone independence which did not satisfy the commutativity axiom. A new version of the axiomatic theory that would include monotone independence was developed by Muraki [18], who showed that when we drop the commutativity axiom, we have five nice notions of independence 
(in addition to the former three notions, there is monotone independence and its twin version, anti-monotone independence).
As concerns their relation to tensor independence, Franz showed in [6] that one can generalize the tensor reduction approach developed in [11,12] to include monotone and anti-monotone independences.
Other notions of independence have also attracted considerable interest 
(conditional freeness, freeness with subordination, orthogonal independence, conditionally monotone independence, matricial freeness, $q$-independence, infinitesimal 
freeness, etc.). 

The first motivation of this paper comes from the fact that each notion of independence leads to a different noncommutative probability theory, which 
does not seem very satisfactory. Therefore, relations between these notions are of interest. Likewise, a possibility to include them in one unified framework should be of importance [11,15]. For this reason, we return to the idea of `reducing' other notions of independence to tensor independence. In this paper, we focus on
the concept of conditionally monotone (shortly, c-monotone) independence of Hasebe [8,9] which generalizes monotone independence. 

As we have already remarked, boolean, free, monotone and anti-monotone independences have been reduced
to tensor independence on suitably constructed larger algebras. 
In fact, freeness was included in this framework in the limit sense [11] and
the actual `reduction' of freeness to tensor independence [12], 
where free random variables 
take the form of infinite series, can also be easily generalized to include
conditional freeness. In turn, in this paper, we show that the conditionally monotone independence 
can also be `reduced' to tensor independence. 

In particular, we demonstrate that if $a_1,a_2$ are two variables which are c-monotone with respect to a pair of states $(\varphi,\psi)$, then variables of the form
$$
A_1=a_1\otimes r_1\;\;\;{\rm and}\;\;\;A_2=s_1\otimes a_2'+s_2\otimes a_2''
$$
have the same distributions as $a_1$ and $a_2$, respectively, w.r.t. to a pair of suitably defined tensor product states and, moreover, are c-monotone with respect to that pair, 
where $a_2',a_2''$ are two tensor independent copies of $a_2$, whereas $r_1,s_1,s_2$ are certain projections, such that $s_1\perp s_2$. 

The second motivation concerns the associated products of graphs. It is interesting to observe that one can associate nice products of rooted graphs with the main notions of noncommutative independence:
\begin{eqnarray*}
tensor\;independence &\rightsquigarrow& Cartesian\;product,\\
freeness&\rightsquigarrow& free\;product,\\
monotone\;independence&\rightsquigarrow&  comb\;product,\\
boolean\;independence&\rightsquigarrow& star\;product,\\
orthogonal\;independence&\rightsquigarrow&  orthogonal \;product,\\
s\!-\!freeness &\rightsquigarrow& s\!-\!free\;product,\\
c\!-\!monotone\;independence &\rightsquigarrow& c\!-\!comb\;product,
\end{eqnarray*}
where the {\it c-comb product of birooted graphs} (graphs with two distinguished vertices called roots) is defined in this paper.

These graph products give a good intuitive understanding of these notions. For instance, it was clear from the beginning of free probability that the free product of graphs 
(especially the Cayley graphs of free products of groups) corresponds to freeness (for the construction of the free product of arbitrary rooted graphs, see [21]). Then, it was shown in [1] that the comb product of rooted graphs corresponds to monotone independence.
Next, it was observed by the author that the star product of rooted graphs corresponds to boolean independence. This observation was used by Obata [19]
to a study of spectral properties of this product.
For more on all these products and their spectral properties, see the monograph of Hora and Obata [10]. Finally, products of rooted graphs related to orthogonal independence and freeness with subordination 
(shortly, s-freeness) were introduced in [13] and further studied in [2,14]. It is noteworthy that all these products are related to subordination and appear in the decompositions of the free product of graphs [2,13].

Let us also remark that the tensor product realizations are canonical for 
the adjacency matrices of these product graphs since they show quite clearly how to glue copies of one graph to copies of another graph, as we showed in 
[13] and [2].
Similarly, the tensor product realization of c-monotone independence leads 
in a natural way to a new product of birooted graphs called the c-comb product of birooted graphs, introduced and studied in this paper. Its pair of spectral distributions is related to the c-comb convolution of probability measures on the real line.  
A suitable modification of the tensor product realization leads to the 
{\it c-comb loop product of birooted graphs} whose pair of spectral distributions is 
related to the multiplicative c-monotone convolution.

Whenever we speak of a graph, we mean a uniformly locally finite rooted (or, birooted) 
graph. In this paper, a graph can be disconnected and it can have loops.

\section{Extensions of algebras and states}

By a noncommutative probability space we understand a 
pair $({\mathcal A}, \varphi)$, where ${\mathcal A}$ is a unital algebra and 
$\varphi$ is a normalized linear functional (i.e. such that $\varphi(1)=1$). In the case when ${\mathcal A}$ is a unital *-algebra, then we also 
require that $\varphi$ is positive, i.e $\varphi(xx^*)\geq 0$ for any $x\in {\mathcal A}$. Then the pair $({\mathcal A}, \varphi)$ is called a *-noncommutative probability space. In both situations we will call 
$\varphi$ a state. If ${\mathcal A}$ is equipped with two states, $\varphi$ and $\psi$, then the triple $({\mathcal A}, \varphi, \psi)$ 
will still be considered to be a (*-) noncommutative probability space.

The `reduction' of boolean, free and monotone independences to tensor independence is done by constructing suitable extensions
of the considered noncommutative probability spaces [11]. 

Let us recall the definitions of boolean, monotone, free and orthogonal independences (see, for instance, [13]).
Let $({\mathcal A}, \varphi)$ be a noncommutative probability space.
\begin{enumerate}
\item
The family $\{{\mathcal A}_{j}: j\in I\}$ of subalgebras of ${\mathcal A}$ is {\it boolean independent} with respect to $\varphi$ if 
$$
\varphi(a_1\cdots a_n)=\varphi(a_1)\cdots \varphi(a_n)
$$
for any $a_i\in \mathcal{A}_{j_i}$, where $j_1\neq \cdots \neq j_n$. 
\item
The family $\{{\mathcal A}_{j}: j\in I\}$ of subalgebras of ${\mathcal A}$, where $I$ is linearly ordered, is 
{\it monotone independent} with respect to $\varphi$ if 
$$
\varphi(\cdots a_{i-1}a_{i}a_{i+1}\cdots )=\varphi(a_i)\varphi (\cdots a_{i-1}a_{i+1}\cdots)
$$ 
whenever $a_i\in \mathcal{A}_{j(i)}$, where $j_1\neq \cdots \neq j_n$ and $j_{i-1}<j_i$, $j_{i}>j_{i+1}$, 
with the understanding that if $i=1$, then $j_i>j_{i+1}$, and if 
$i=n$, then $j_{i-1}<j_i$.
\item
The family $\{{\mathcal A}_{j}: j\in I\}$ of unital subalgebras of ${\mathcal A}$ is
{\it free} with respect to $\varphi$ if
$$
\varphi(a_1\cdots a_n)=0
$$
whenever $a_i\in {\mathcal A}_{j_i}\cap {\rm Ker}\varphi$ for any $i$.
\item
Let $\psi$ be another linear normalized functional on ${\mathcal A}$ and let ${\mathcal A}_{1}, {\mathcal A}_{2}$ be subalgebras
of ${\mathcal A}$. Then ${\mathcal A}_{2}$ is {\it orthogonal} to ${\mathcal A}_{1}$ with respect to $(\varphi, \psi)$ if
$$
\varphi(\cdots a_{i-1}a_{i}a_{i+1}\cdots )=
\psi(a_{i})\left(\varphi(\cdots a_{i-1}a_{i+1}\cdots )- \varphi(\cdots a_{i-1})\varphi(a_{i+1}\cdots )\right)
$$
for any $a_{i-1},a_{i+1}\in {\mathcal A}_{1}$, $a_{i}\in {\mathcal A}_{2}$, and
$\varphi(a_{1}\cdots)=\varphi(\cdots a_n)=0$ whenever $a_1,a_n\in {\mathcal A}_{2}$.
\end{enumerate}

\begin{Definition}
{\rm 
Let $({\mathcal A}, \varphi)$ be a noncommutative probability space. 
Define its extension $(\widetilde{\mathcal A}, \widetilde{\varphi})$, where $\widetilde{\mathcal A}={\mathcal A}(P)$, the (unital) algebra generated by ${\mathcal A}$ and the indeterminate $P$, such that $1P=P1=P^2=P$, where $1$ is 
the unit in ${\mathcal A}$ and $\widetilde{\varphi}:\widetilde{{\mathcal A}}\rightarrow {\mathbb C}$ is the linear functional given by the linear extension of $\widetilde{\varphi}(P)=1$ and 
$$
\widetilde{\varphi}(P^{\alpha}a_1Pa_2P\cdots Pa_{n}P^{\beta})=\varphi(a_1)\varphi(a_2)\cdots \varphi(a_n)
$$
where $\alpha, \beta\in \{0,1\}$ and $a_1, \ldots, a_n\in {\mathcal A}$. It is easy to see that $\widetilde{\varphi}$ is a normalized linear functional on $\widetilde{\mathcal A}$ and that
$\widetilde{\varphi}(cPc')=\widetilde{\varphi}(c)\widetilde{\varphi}(c')$ for any $c,c'\in \widetilde{\mathcal A}$.
}
\end{Definition}

If ${\mathcal A}$ is a *-algebra, we assume in addition that $P^*=P$. 
It is known that in that case $\widetilde{\varphi}$ is positive if $\varphi$ is positive and hence a state on 
$\widetilde{\mathcal A}$ [11]. If ${\mathcal A}$ is a $C^*$-algebra and 
$({\mathcal H}, \pi, \xi)$ is the GNS triple 
associated with $(\mathcal{A}, \varphi)$, then $({\mathcal H}, \widetilde{\pi}, \xi)$ is the GNS triple
associated with $(\widetilde{\mathcal{A}}, \widetilde{\varphi})$, where $\widetilde{\pi}(a)=\pi(a)$ 
for any $a\in {\mathcal A}$ and $\widetilde{\pi}(P)=P_{\xi}$, 
the orthogonal projection onto ${\mathbb C}\xi$. Note that if we have two different states on ${\mathcal A}$, 
say $\varphi$ and $\psi$, then we have two different GNS triples, say $({\mathcal H}_{1}, \pi_1, \xi_2)$ and 
$({\mathcal H}_{2}, \pi_2, \xi_2)$, respectively. It is possible to 
take ${\mathcal H}_{1}={\mathcal H}_{2}$ and $\xi_1=\xi_2$ 
(by taking the tensor product of the original Hilbert spaces) at the expense of losing cyclicity of the representations.

\begin{Proposition}
Let $({\mathcal A}_{1},\varphi_1), ({\mathcal A}_{2}, \varphi_2)$ be noncommutative probability spaces
and let $\widetilde{{\mathcal A}}_{1}={\mathcal A}(P)$, $\widetilde{{\mathcal A}_{2}}={\mathcal A}(Q)$. 
Let $a_1\in {\mathcal A}_{1}$, $a_2\in {\mathcal A}_{2}$ and 
$\varphi=\widetilde{\varphi}_1\otimes \widetilde{\varphi}_2$.
\begin{enumerate}
\item
The variables
$$
A_1=a_1\otimes Q\;\;\;{\rm and}\;\;\;
A_2=P\otimes a_2
$$
are boolean independent with respect to $\varphi$.
\item
The variables
$$
A_1=a_1\otimes Q\;\;\;{\rm and}\;\;\;
A_2=1_1\otimes a_2
$$
where $1_1$ is the unit in ${\mathcal A}_{1}$, are monotone independent with respect to $\varphi$.
\item
The variables
$$
A_1=a_1\otimes Q\;\;\;{\rm and}\;\;\;
A_2=P^{\perp}\otimes a_2,
$$
where $P^{\perp}=1-P$, are orthogonally independent with respect to $\varphi$. 
\end{enumerate}
\end{Proposition}
{\it Proof.}
For the proofs of (1),(2) and (3) we refer the reader to [11], [6] and [13], respectively.
\hfill $\blacksquare$\\

Let us observe that in the case of monotone independence of 
two algebras, it suffices to extend $({\mathcal A}_{2},\varphi_2)$, but in the case of a general 
family of noncommutative probability spaces $\{({\mathcal A}_{j}, \varphi_j):j\in J\}$ one needs to extend
each $(\mathcal{A}_{j}, \varphi_j)$. Let us add that the above tensor realizations can be generalized to 
arbitrary families of noncommutative probability spaces (in the monotone case, they have to be indexed by a 
linearly ordered set since only for such families monotone independence is meaningful). For details, see [11] and [6].

A much more sophisticated realization was found for free random variables since one has to take infinite sequences of copies of $a_{1}\in {\mathcal A}_{1}$ 
and $a_{2}\in {\mathcal A}_{2}$,  
namely $(a_{1,n})_{n\geq 1}$ and $(a_{2,n})_{n\geq 1}$. Then
the variables
$$
A_1=\sum_{n=1}^{\infty}a_{1,n}\overline{\otimes} q_n\;\;\;{\rm and}\;\;\;A_2=\sum_{n=1}^{\infty}p_n\overline{\otimes} a_{2,n},
$$
where $\{p_n:n\in \mathbb{N}\}$ and $\{q_n: n\in \mathbb{N}\}$ are suitably defined 
sequences of orthogonal projections that give (strongly convergent) resolutions of the identity, 
are free with respect to the tensor product state 
$\varphi=\widetilde{\varphi}_1^{\otimes \infty}\otimes \widetilde{\varphi}_2^{\otimes \infty}$
Here, $\overline{\otimes}$ is a suitably defined tensor product (we used Berberian's theory to define a suitable notion of
convergence in *-algebras in which $\overline{\otimes}$ reminds the von-Neumann tensor product).

\section{C-monotone independence}
Let us recall the definition of c-monotone independence of Hasebe [8]. Our definition is simpler in the sense that we reduce Hasebe's definition to a `local maximum' condition for indices with the tacit understanding that this `local maximum' can also
occur at the beginning or at the end of the tuple of indices.
\begin{Definition}
{\rm
Let $(A, \varphi, \psi)$ be a unital algebra equipped with two states. Let $J$ be a linearly ordered set. 
A family of subalgebras $\{A_{j}:j\in J\}$  is said to be {\it c-monotone independent}  
with respect to the pair $(\varphi, \psi)$ if and only if 
it holds that
\begin{eqnarray*}
\varphi (a_1\cdots a_n)&=&
(\varphi(a_i)-\psi(a_i)) \varphi(a_1\cdots a_{i-1})\varphi(a_{i+1}\cdots a_{n})\\
&& +
\psi(a_{i})\varphi(a_1\cdots a_{i-1}a_{i+1}\cdots a_n)
\end{eqnarray*}
for any $a_k\in {\mathcal A}_{j_k}$, where $j_1, \ldots, j_n\in J$ and $j_{i-1}<j_i>j_{i+1}$, and we understand that
if $i=1$ or $i=n$, then only the right or left inequality is required to hold, respectively (then $\varphi$ evaluated at the empty product of elements 
is set to be equal to one) and that these subalgebras are monotone independent with respect to $\psi$.
}
\end{Definition}

As compared with Hasebe's definition [8, Definition 3.8], instead of three conditions: for $j=1$, $1<j<n$ and $j=n$, we give one `local maximum' condition, with the understanding that 
if $j=1$ or $j=n$, then the moments of empty products of variables 
that precede $a_1$ or follow $a_n$, respectively, are equal to one.

We would like to reduce c-monotone independence to tensor independence. Before we do it in all generality, let us consider the case of two algebras equipped with two states: $({\mathcal A}_{1},\varphi_1, \psi_1)$ and $(\mathcal{A}_{2},\varphi_2,\psi_2)$, and let 
$\widetilde{\mathcal A}_{1}={\mathcal A}_{1}(P)$ and $\widetilde{\mathcal A}_{2}={\mathcal A}_{2}(Q)$,
with extended states $\widetilde{\varphi}_{1},\widetilde{\varphi}_{2}, 
\widetilde{\psi}_{1}, \widetilde{\psi}_{2}$. 
In this setting, consider the tensor product states
\begin{eqnarray*}
\varphi&=&\widetilde{\varphi}_{1}\otimes 
\widetilde{\varphi}_2\otimes \widetilde{\psi}_2\\
\psi&=&\widetilde{\psi}_{1}\otimes 
\widetilde{\psi}_2\otimes \widetilde{\psi}_2
\end{eqnarray*}
on the tensor product algebra 
${\mathcal A}:=\widetilde{\mathcal A}_{1}\otimes \widetilde{\mathcal A}_{2}\otimes \widetilde{\mathcal A}_{2}$.

\begin{Proposition} Let $j_i:{\mathcal A}_{i}\rightarrow {\mathcal A}$ 
be homomorphisms given by $j_i(a_i)=A_i$, where $i=1,2$ and we set
\begin{eqnarray*}
A_1&=&a_1\otimes Q\otimes Q\\
A_2&=&P\otimes a_2 \otimes 1_2 +P^{\perp}\otimes 1_2\otimes a_2
\end{eqnarray*}
for $a_i\in {\mathcal A}_i$. 
Then $\varphi \circ j_i=\varphi_i$ and $\psi\circ j_i=\psi_i$ for $i=1,2$
and the subalgebras of ${\mathcal A}$ generated by $A_1$ and $A_2$, respectively, 
are c-monotone independent with respect to $(\varphi, \psi)$. 
\end{Proposition}
{\it Proof.}
Note that $j_1$ is a non-unital homomorphism which maps $1_{1}$ onto $1_{1}\otimes Q\otimes Q$, whereas $j_2$ is a unital homomorphism which maps $1_{2}$ onto
$1_{1}\otimes 1_{2}\otimes 1_{2}$.
Using the definition of the extended states, it is easy to see that 
${\varphi}(A_i)=\varphi_{i}(a_i)$ and ${\psi}(A_i)=\psi_{i}(a_i)$, where $i=1,2$. Therefore,
$\varphi \circ j_i=\varphi_i$ and $\psi\circ j_i=\psi_i$.
Now, let $A=j_1(a), A'=j_1(a'), B=j_2(b)$, where $a,a'\in {\mathcal A}_{1}$ and 
$b\in {\mathcal A}_{2}$, and let $C,C'$ be simple tensors from ${\mathcal A}$ such that $C$ ends with an element from ${\mathcal A}_{2}$ and $C'$ begins with an element from ${\mathcal A}_2$. In general, they are linear combinations of simple tensors, but in our computations it suffices to take simple tensors. 
We have
\begin{eqnarray*}
\varphi(CABA'C')&=&
\varphi(C(a\otimes Q\otimes Q)(P\otimes b\otimes 1_2
+P^{\perp}  \otimes 1_2 \otimes b)
(a'\otimes Q\otimes Q)C')\\
&=&
\varphi(C(aPa'\otimes QbQ \otimes Q)C')+\varphi(C(aP^{\perp}a'\otimes Q\otimes QbQ)C')\\
&=& 
{\varphi}_{2}(b)\varphi(C (aPa'\otimes Q\otimes Q)C')+\psi_2(b)\varphi(C (aa'\otimes Q\otimes Q)C')\\
&&-\psi_2(b)\varphi(C (aPa'\otimes Q\otimes Q)C')\\
&=&
({\varphi}(B)-\psi(B))\varphi(CA)\varphi(A'C')+
\psi(B)\varphi(C AA'C')
\end{eqnarray*}
where we used the definition of the extended states repeatedly. In particular, in the last equation 
we used
$$
\varphi(D(P\otimes Q\otimes Q)D')=\varphi(D)\varphi(D')
$$
for any $D,D'\in {\mathcal A}$ and ${\varphi}(B)=\varphi_{2}(b)$, ${\psi}(B)=\psi_{2}(b)$.
This gives the desired expression of Definition 3.1 in the case when the 'local maximum' is in the middle of the moment, 
namely $1<i<n$. When $i=1$ or $i=n$, computations are easy: for instance, 
if $i=1$, there is no $A$ or $C$ and the term with $P^{\perp}$ can be deleted since $\widetilde{\varphi}_1(P^{\perp}a'c')=0$.
Finally, it is easy to see that when we take $\psi$ instead of $\varphi$, then we obtain $\varphi(B)=\psi(B)$ in the above computation, which gives the defining condition of monotone independence. This completes the proof. \hfill $\blacksquare$

\begin{Remark}
{\rm Let us make some observations concerning c-monotone independence and its tensor realization.
\begin{enumerate}
\item
In the case of two algebras, it sufffices to use tensors of order three (one tensor for the first algebra and two tensors for the second one), as we did in Proposition 3.1. However, as will be seen below, in the case of a general totally ordered index set $J$, this assymetry will not be visible since 
usually each index $j\in  J$ is smaller than some indices and bigger than some other ones. Of course, if $J=\{1,2\}$, we can also use states of the form
$$
\varphi=\widetilde{\varphi}_{1}\otimes\widetilde{\psi}_{1}\otimes \widetilde{\varphi}_{2}\otimes \widetilde{\psi}_{2},
\,\,\,
\psi=\widetilde{\psi}_{1}\otimes\widetilde{\psi}_{1}\otimes \widetilde{\psi}_{2}\otimes \widetilde{\psi}_{2},
$$ 
but we shall obtain the same mixed moments w.r.t. $\varphi$ and $\psi$ as in the case of tensor products of order three.
\item
It follows from the proof of Proposition 3.1 that the same mixed moments with respect to 
$\varphi$ and $\psi$ are obtained when we take
\begin{eqnarray*}
A_2&=&P\otimes a_2 \otimes Q +P^{\perp}\otimes Q\otimes a_2
\end{eqnarray*}
instead of the $A_2$ in Proposition 3.1.
This form of $A_2$ is slightly more convenient in the construction of the associated products of graphs.
\item
Informally, c-monotone independence is a mixture of boolean independence and orthogonal independence. Namely, by Proposition 3.1, the variable $A_2$ is a sum of two (tensor independent) copies of $a_2$, 
one of which forms with $A_1$ a boolean independent pair w.r.t. $\varphi$, whereas the second one forms an orthogonally independent pair w.w.t. $\varphi$.
This aspect will also appear in Section 3, where we discuss the c-monotone product of rooted graphs.
\end{enumerate}
}
\end{Remark}

\begin{Example}
{\rm 
The lowest order nontrivial mixed moment is of the form:
\begin{eqnarray*}
\varphi(ABA')
&=&
(\widetilde{\varphi}_{1}\otimes 
\widetilde{\varphi}_2\otimes \widetilde{\psi}_2 )
((a\otimes Q\otimes Q)(P\otimes b\otimes 1_2\\
&&+
P^{\perp}\otimes 1_2 \otimes b)(a'\otimes Q \otimes Q)) \\
&=&
\widetilde{\varphi}_{1}(aPa')\widetilde{\varphi}_{2}(QbQ)+
\widetilde{\varphi}_{1}(aP^{\perp}a')\widetilde{\psi}_{2}(QbQ)\\
&=&
{\varphi}_{1}(a)\varphi_{1}(a')\varphi_{2}(b)+
{\varphi}_{1}(aa')\psi_2(b)-
{\varphi}_{1}(a)\varphi_1(a')\psi_{2}(b).
\end{eqnarray*}
It can be easily seen that this expression agrees with the mixed moment 
of $aba'$ if $a,a'\in {\mathcal A}_{1}$ and ${\mathcal A}_{2}$ and we use Definition 3.1.
}
\end{Example}

\begin{Example}
{\rm 
Let us also compute a mixed moment $\varphi(ABA'B'A'')$. Observe that a nonvanishing contribution is given by the following 
product:
\begin{eqnarray*}
&&(a\otimes Q\otimes Q)(P\otimes b\otimes 1_2+
P^{\perp}\otimes 1_2 \otimes b)\\
&&\times 
(a'\otimes Q \otimes Q)
(P\otimes b'\otimes 1_2+
P^{\perp}\otimes 1_2 \otimes b')
(a''\otimes Q \otimes Q) \\
&=&
aPa'Pa''\otimes QbQb'Q\otimes Q+
aPa'P^{\perp}a''\otimes QbQ\otimes Qb'Q\\
&&+
aP^{\perp}a'Pa''\otimes Qb'Q\otimes QbQ+
aP^{\perp}a'P^{\perp}a''\otimes Q\otimes QbQb'Q
\end{eqnarray*}
which gives
\begin{eqnarray*}
&&\varphi_{1}(a)\varphi_1(a')\varphi_1(a'')\varphi_2(b)\varphi_2(b')
+\varphi_{1}(a)(\varphi_1(a'a'')-\varphi_1(a')\varphi_1(a''))\varphi_2(b)\psi_2(b')\\
&&+(\varphi_{1}(aa')-\varphi_1(a)\varphi_1(a'))\varphi_1(a'')\varphi_2(b')\psi_2(b)\\
&&+(\varphi_{1}(aa'a'')-\varphi_1(aa')\varphi_1(a'')-\varphi_1(a)\varphi_1(a'a'')+\varphi_1(a)\varphi_1(a')\varphi_1(a''))
\psi_2(b)\psi_2(b').
\end{eqnarray*}
Again, it can be verified, using Definition 3.1, that this moment agrees with the moment of conditionally monotone independent random variables.
}
\end{Example}

We proceed with proving the theorem in the general case. It is less explicit than that of Proposition 3.1,
but is is relatively simple. We will use the tensor product algebra of the form
$$
{\mathcal A}=\bigotimes_{j\in J}(\widetilde{\mathcal A}_{j}\otimes \widetilde{\mathcal A}_{j})
$$ 
and the associated tensor product states
\begin{eqnarray*}
\varphi&=&\bigotimes_{j\in J}(\widetilde{\varphi}_{j}\otimes \widetilde{\psi}_{j})\\
\psi&=&\bigotimes_{j\in J}(\widetilde{\psi}_{j}\otimes \widetilde{\psi}_{j}).
\end{eqnarray*}
Depending on $j$, we will use decompositions of ${\mathcal A}$ of the form
$$
{\mathcal A}=(\widetilde{\mathcal A}_{j}\otimes \widetilde{\mathcal A}_{j})\otimes\bigotimes_{k\neq j}(\widetilde{\mathcal A}_{k}\otimes 
\widetilde{\mathcal A}_{k})
$$
for any $j\in J$, depending on to which algebra we want to associate an element of the tensor product.
Therefore, the order in which tensor components of a given vector are written is not fixed. 

\begin{Remark}
{\rm 
Let us observe that the projection of the form
$$
P:=\bigotimes_{j\in J} (P_j\otimes P_j)
$$
separates all elements of the algebra ${\mathcal A}$ in the sense that 
$$
\varphi(A_{j_1}\cdots A_{j_{k}}PA_{j_{k+1}}\cdots A_n)=\varphi(A_{j_1}\cdots A_{j_{k}})\varphi(A_{j_{k+1}}\cdots A_n)
$$
since $P_j$ separates all elements of the algebra ${\mathcal A}_{j}$. Moreover,
$$
\widetilde{\varphi}_{i}(a'P_jaP_ja'')=\varphi_{j}(a)\widetilde{\varphi}_{j}(a'P_{j}a'')
$$
for any $a,a',a''\in {\mathcal A}_{j}$ and any $j\in J$.
Therefore, in particular, whenever the product $A_{j_1}\cdots A_{j_n}$ contains a product 
$P_{j}a_{j}P_{j}$ for some $j=j_k$ at some tensor site of color 1 (odd) or color 2 (even), we can replace it by 
$$
\varphi_{j}(a_{j})P_{j}=\varphi(A_{j})P_j \;\;\;{\rm or} \;\;\;\psi_{j}(a_{j})P_j=\psi(A_{j})P_j,
$$
respectively.
}
\end{Remark}

\begin{Theorem}
The subalgebras of ${\mathcal A}=\bigotimes_{j\in J}(\widetilde{\mathcal A}_{j}\otimes \widetilde{\mathcal A}_{j})$ 
generated by elements of the form
\begin{eqnarray*}
A_j&=&(a_j\otimes 1_j)\otimes p_{j}+(1_j\otimes a_j)\otimes (p_{j}'-p_{j}),
\end{eqnarray*}
where  $a_j\in {\mathcal A}_j$ and $j\in J$, respectively, and 
\begin{eqnarray*}
p_j&=&\bigotimes_{k\neq j}(P_{k}\otimes P_{k})\\
p_j'&=& \bigotimes_{k<j}(1_{k}\otimes 1_k)\otimes \bigotimes_{k>j}(P_{k}\otimes P_k),
\end{eqnarray*}
are c-monotone independent with respect to $(\varphi, \psi)$.
\end{Theorem}
{\it Proof.}
Let $a\in {\mathcal A}_{i}, b\in {\mathcal A}_{k}, c\in {\mathcal A}_{l}$, where 
$i<j$ and $j>k$. Consider the associated elements of ${\mathcal A}$:
\begin{eqnarray*}
A&=&(a\otimes 1_{i})\otimes p_i + (1_i\otimes a)\otimes (p_i'-p_i)\\
B&=&(b\otimes 1_{j})\otimes p_j + (1_j\otimes b)\otimes (p_j'-p_j)\\
C&=&(c\otimes 1_{k})\otimes p_k + (1_k\otimes c)\otimes (p_k'-p_k)
\end{eqnarray*} 
and write $B=B'+B''$, where  
\begin{eqnarray*}
B'&=&(b\otimes 1_{j})\otimes p_j\\
B''&=&(1_j\otimes b)\otimes (p_j'-p_j).
\end{eqnarray*}
Using Remark 3.2, we obtain
$$
AB'C=A\left((P_jbP_j\otimes P_j)\otimes p_j\right)C=\varphi(B)APC
$$
due to the fact that $i<j>k$ and thus elements $A,C$ contain projections at all sites 
associated with $j$.
Similarly, we have
\begin{eqnarray*}
AB''C&=&A\left((1_j\otimes b)\otimes (p_j'-p_j)\right)C\\
&=& A\left((P_j\otimes P_jbP_j)\otimes p_j'\right)C-A\left((P_j\otimes P_jbP_j)\otimes p_j\right)C\\
&=&\psi(B)A\left((1_j\otimes 1_j)\otimes p_j'\right)C-\psi(B)APC.
\end{eqnarray*}
Now, $p_j'$ has projections at all sites associated with indices $l>j$, but at these sites both $A$ and $C$ have
projections as well, hence
$$
A(1_j\otimes 1_j\otimes p_j')C=AC
$$
and thus
$$
ABC=(\varphi(B)-\psi(B))APC +\psi(B)AC
$$
which gives the required recurrence for moments since
\begin{eqnarray*}
\varphi(\cdots ABC \cdots )&=&
(\varphi(B)-\psi(B))\varphi(\cdots APC\cdots ) +\psi(B)\varphi(\cdots AC\cdots)\\
&=&(\varphi(B)-\psi(B))\varphi(\cdots A)\varphi(C\cdots ) +\psi(B)\varphi(\cdots AC\cdots)
\end{eqnarray*}
and that completes the proof when $B$ is in the middle of the moment. The case when
$B$ is the first or the last element in the product reduces to the above when we 
set $A=1$ or $C=1$, respectively. Therefore, the proof is completed.
\hfill $\blacksquare$

\section{C-comb product of graphs}

Using c-monotone independence and its tensor product realization, 
we can define a new product of rooted graphs, called the c-comb product of rooted graphs.
It is a generalization of the comb product of rooted graphs studied in [1,10] and it is related to the orthogonal product of rooted graphs, which turned out important in the decompositions of the free product of graphs [2,13].

\begin{Remark}
{\rm 
Let us recall some definitions and facts on products of rooted graphs.
\begin{enumerate}
\item
By a {\it rooted graph} we understand a pair $({\mathcal G},e)$, where
${\mathcal G}=(V,E)$ is a non-oriented graph with the set of vertices $V$ and the set of edges $E$, and $e\in V$ is a distinguished vertex called the {\it root}. We assume that our rooted graphs are uniformly locally finite. They may be disconnected, but
we will use only graphs with a finite number of components (the root will belong 
to one of these components). To each vertex $x\in V$ we associate a unit vector $\delta(x)$ in the 
Hilbert space $\ell^{2}(V)$.
The orthogonal projection onto ${\mathbb C}\delta(x)$ will be denoted $P_{x}$ and
$P_{x}^{\perp}=1-P_{x}$. 
\item
If $(({\mathcal G}_{1}, e_1), ({\mathcal G}_{2}, e_{2}))$ is an ordered pair
of rooted graphs, then by their {\it disjoint union} we 
will understand the graph with two distinguished vertices
$$
({\mathcal G}_{1}\oplus {\mathcal G}_{2},e_1,e_2)\equiv
({\mathcal G}_{1},e_1)\oplus ({\mathcal G}_{2},e_2)
$$
with the set of vertices $V=V_1\cup V_2$, the set of edges $E=E_{1}\cup E_{2}$ and the root chosen to be $e=e_1$. For simplicity, we will also write 
${\mathcal G}_{1}\oplus {\mathcal G}_{2}$. In view of the root selection, the order in the disjoint union of {\it rooted} graphs is important. If $a_1=a({\mathcal G}_{1})$ 
and $a_2=a({\mathcal G}_{2})$ are the adjacency matrices of ${\mathcal G}_{1}$ and ${\mathcal G}_{2}$, respectively, then the matrix of the form
$$
a=\left(\!\begin{array}{cc}
a_1 & 0\\
0 & a_2
\end{array}
\!\right)
$$
is the adjacency matrix of ${\mathcal G}_{1}\oplus {\mathcal G}_{2}$ and we will write $a=a_{1}\oplus a_{2}$.
\item
The {\it comb product} of rooted graphs $({\mathcal  G}_{1},e_1)$ and $({\mathcal  G}_{2},e_2)$
is the rooted graph
obtained by attaching a copy of ${\mathcal  G}_{2}$ at its root $e_{2}$
to each vertex of ${\mathcal  G}_{1}$, where we denote by $e$ the vertex obtained
by identifying $e_1$ and $e_2$. We will denote it by
$({\mathcal  G}_{1},e_1)\vartriangleright ({\mathcal  G}_{2},e_2)=
({\mathcal  G}_{1}\vartriangleright {\mathcal  G}_{2},e)$, or simply 
${\mathcal  G}_{1}\vartriangleright {\mathcal  G}_{2}$, if it is clear which vertices are 
the roots.
If we identify the set of vertices 
with $V_{1}\times V_{2}$, then $e$ is identified with $e_{1}\times e_{2}$.
If $({\mathcal  G}_{1},e_{1})$ and $({\mathcal  G}_{2},e_{2})$ have 
adjacency matrices $a_1$ and $a_1$ and spectral 
distributions $\mu_1$ and $\mu_2$, respectively (these are distributions defined by the moments of these matrices w.r.t. the states defined by vectors 
$\delta(e_1)$ and $\delta(e_2)$, respectively),
then the adjacency matrix of their comb product is of the form
$$
a({\mathcal  G}_{1}\vartriangleright {\mathcal  G}_{2}) =S_1+S_2
$$
where 
$$
S_1=a_{1}\otimes P_{e_2}\;\;\;{\rm and}\;\;\;S_2=1\otimes a_{2},
$$
living in $\ell^2(V_1)\otimes \ell^2(V_2)$,
and the spectral distribution of ${\mathcal  G}_{1}\vartriangleright {\mathcal  G}_{2}$ 
is given by the monotone additive convolution $\mu_1 \vartriangleright \mu_2$ of Muraki [17].
Recall that 
$$
F_{\mu_1\vartriangleright \mu_2}(z)=F_{\mu_1}(F_{\mu_2}(z)),
$$
where
$F_{\mu}(z)=1/G_{\mu}(z)$ is the reciprocal Cauchy transform of $\mu$. 
For more details on the monotone additive convolution in the context of the comb product of rooted graphs, see [1,10].
\item
The {\it star product} of rooted graphs $({\mathcal  G}_{1},e_1)$ and $({\mathcal  G}_{2},e_2)$
is the rooted graph obtained by attaching ${\mathcal  G}_{2}$ at its root 
to the root of ${\mathcal  G}_{1}$. We denote it by $({\mathcal  G}_{1},e_1)\star ({\mathcal  G}_{2},e_2)=
({\mathcal  G}_{1}\star {\mathcal  G}_{2},e)$,
where we denote by $e$ the vertex obtained by identifying $e_1$ and $e_2$, or simply 
${\mathcal  G}_{1}\star {\mathcal  G}_{2}$. If we identify its set of vertices $V$
with $V_1\star V_2:=(\{e_1\}\times V_{2})\cup (V_{1}\times \{e_2\})$, then $e$ is identified with $(e_{1},e_{2})$. Its
adjacency matrix is of the form
$$
a({\mathcal  G}_{1}\star {\mathcal  G}_{2}) =a_1\otimes P_{e_{2}}+P_{e_1}\otimes a_{2}
$$
living in $\ell^2(V)\subset \ell^2(V_1)\otimes \ell^2(V_2)$,
and the spectral distribution of ${\mathcal  G}_{1}\star {\mathcal  G}_{2}$ is given by
$$
F_{\mu_1\uplus \mu_2}(z)=F_{\mu_1}(z)+F_{\mu_2}(z)-z.
$$
\item
The {\it orthogonal product} of rooted graphs $({\mathcal  G}_{1},e_{1})$ and $({\mathcal  G}_{2},e_{2})$ 
is the rooted graph obtained by attaching a copy of ${\mathcal  G}_{2}$ at
its root $e_{2}$ to each vertex of ${\mathcal  G}_{1}$ but the root $e_{1}$, with the root of the product taken to be $e=e_{1}$.
We denote it by $({\mathcal  G}_{1},e_1)\vdash ({\mathcal  G}_{2},e_2)=({\mathcal G}_{1}\vdash {\mathcal G}_{2},e)$, 
or simply ${\mathcal  G}_{1}\vdash {\mathcal  G}_{2}$ if it is clear which vertices are the roots.
If its set of vertices $V$ is identified with 
$V_1\vdash V_2:=(V_{1}^{0}\times V_{2})\cup \{(e_{1},e_{2})\}$, where
$V_{1}^{\circ}=V_1\setminus \{e_1\}$, then $e$ is identified with 
$(e_{1},e_{2})$. 
The adjacency matrix of $({\mathcal  G}_{1}\vdash {\mathcal  G}_{2},e)$ is of the form
$$
a({\mathcal  G}_{1}\vdash {\mathcal  G}_{2}) =a_1\otimes P_{e_{2}}+P_{e_1}^{\perp}\otimes a_{2}
$$
living in $\ell^2(V)\subset\ell^{2}(V_1)\otimes \ell^2(V_2)$, and the spectral distribution of ${\mathcal  G}_{1}\vdash {\mathcal  G}_{2}$ 
is given by the orthogonal convolution $\mu_1 \vdash \mu_2$ introduced in [13].
It holds that
$$
F_{\mu_1\vdash\mu_2}(z)=F_{\mu_1}(F_{\mu_2}(z))-F_{\mu_2}(z)+z.
$$
For more details on the orthogonal product of rooted graphs, 
see [13] and [2].
\item
Note that both ${\mathcal  G}_{1}\star {\mathcal  G}_{2}$ and ${\mathcal  G}_{1}\vdash {\mathcal  G}_{2}$
are subgraphs of ${\mathcal  G}_{1}\vartriangleright {\mathcal  G}_{2}$. In fact, one can view 
the latter as a `superposition' of the former two products. More formally, 
we have the decomposition
$$
({\mathcal  G}_{1},e_1)\vartriangleright ({\mathcal  G}_{2},e_2)=(({\mathcal G}_{1},e_1)\vdash ({\mathcal G}_{2},e_2))\star ({\mathcal G}_{2},e_2)
$$
and thus it is convenient to treat these adjacency matrices as operators living in 
$\ell^{2}(V_1)\otimes \ell^{2}(V_2)$.
\end{enumerate}
}
\end{Remark}

The above decomposition of the comb product of rooted graphs 
can be generalized in a natural way. Namely, we can distinguish two different vertices in a given graph, say $e,f$,
consider the triple $({\mathcal G}, e,f)$ and call it a {\it birooted graph}, where both $e$ and $f$ are called 
the {\it roots}. By the {\it spectral distribution} of $({\mathcal G}, e,f)$
we will understand the pair $(\mu, \nu)$, where $\mu$ and $\nu$ are spectral distributions of $({\mathcal G},e)$ and $({\mathcal G},f)$, respectively.

\begin{Definition}
{\rm  Let $({\mathcal G}_{1},e_1, f_1)$ and $({\mathcal G}_{2},e_2, f_2)$ 
be two birooted graphs. Let us define a product of rooted graphs $({\mathcal G}_{1},e_1)$ and $({\mathcal G}_{2},e_2)$ 
by the formula
$$
({\mathcal  G}_{1},e_1)\vartriangleright_{f_2}({\mathcal  G}_{2},e_2)=
(({\mathcal G}_{1},e_1)\vdash ({\mathcal G}_{2},f_2))\star ({\mathcal G}_{2},e_2).
$$
In other words, it is the graph obtained by attaching a copy of 
${\mathcal  G}_{2}$ at $e_{2}$ to the root of ${\mathcal  G}_{1}$ and a 
copy of ${\mathcal G}_{2}$ at $f_2$ to the remaining vertices of 
${\mathcal G}_{1}$. If we identify the set of vertices of this graph with 
$V_1\vartriangleright_{f_2} V_2:=(V_1\vdash V_2)\times \{e_2\}\cup \{(e_1,f_2)\}\times V_{2}$, then its root $e$ 
is identified with $(e_1,f_2,e_2)$. We will also use the abbreviated notation $({\mathcal G}_{1}\vartriangleright_{f_2} {\mathcal G}_{2},e)$.
}
\end{Definition}

\begin{Remark}
{\rm Clearly, 
$({\mathcal G}_{1},e_1)\vartriangleright_{e_2} ({\mathcal G}_{2},e_2)=({\mathcal G}_{1},e_1)\vartriangleright ({\mathcal G}_{2},e_2)$, so it is the case when $f_2\neq e_2$ which is of main interest. The main point of the above product is that it is the main component of the c-comb product of graphs, whose definition is based
on the definition of the c-monotone product of states of Hasebe [8]. An example is given in Fig.~1. }
\end{Remark}

\begin{Definition}
{\rm Let $({\mathcal  G}_{1},e_1,f_1)$ and $({\mathcal  G}_{2},e_2,f_2)$ be birooted graphs. By their {\it c-comb product}, denoted 
$({\mathcal G}_{1},e_1,f_1)\vartriangleright_{f_2}
({\mathcal G}_{2},e_2,f_2)$, we understand the disjoint union
$$
({\mathcal G},e,f)=\left(({\mathcal G}_{1},e_1)\vartriangleright_{f_2}({\mathcal G}_{2},e_2)\right)\oplus \left(({\mathcal G}_{1},f_1)\vartriangleright ({\mathcal G}_2,f_2)\right),
$$
where $({\mathcal  G}_{1},e_1)\vartriangleright_{f_2} ({\mathcal  G}_{2},e_2)$
is called the {\it essential component} of $({\mathcal G},e,f)$.
If we identify the set of vertices with 
$V=(V_1\vartriangleright_{f_2} V_2 )\cup (V_1\times V_2)\subseteq (V_{1}\times V_2\times V_2)\cup (V_1\times V_2)$, then
$e$ is identified with $(e_1,f_2,e_2)\in V_{1}\times V_2\times V_2$ 
and $f$ is identified with $(f_1, f_2)\in V_{1}\times V_2$.
}
\end{Definition}

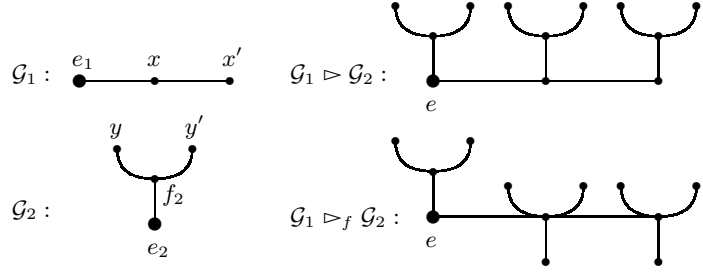
\begin{figure}
\unitlength=1mm
\special{em.linewidth 2pt}
\linethickness{0.5pt}
\begin{picture}(80.00,40.00)(20.00,-8.00)
\put(23.00,18.00){\line(1,0){20.00}}
\put(23.00,18.00){\circle*{1.60}}
\put(33.00,18.00){\circle*{1.00}}
\put(43.00,18.00){\circle*{1.00}}
\put(22.00,20.00){\scriptsize $e_{1}$}
\put(32.00,20.00){\scriptsize $x$}
\put(42.00,20.00){\scriptsize $x'$}
\put(14.00,18.00){\scriptsize ${\mathcal G}_{1}:$}

\put(33.00,-1.00){\circle*{1.60}}
\put(33.00,-1.00){\line(0,1){6.00}}
\put(33.00,5.00){\circle*{1.00}}
\put(28.00,9.00){\circle*{1.00}}
\put(38.00,9.00){\circle*{1.00}}

\put(32.00,-5.00){\scriptsize $e_{2}$}
\put(34.00,2.00){\scriptsize $f_2$}
\put(37.00,11.00){\scriptsize $y'$}
\put(27.00,11.00){\scriptsize $y$}
\put(14.00,0.00){\scriptsize ${\mathcal G}_{2}:$}

\qbezier(28,9)(28,5)(33,5)
\qbezier(33,5)(38,5)(38,9)

\put(69.00,14.00){\scriptsize $e$}
\put(51.00,18.00){\scriptsize ${\mathcal G}_{1}\vartriangleright {\mathcal G}_{2}:$}
\put(70.00,18.00){\line(1,0){30.00}}
\put(70.00,18.00){\line(0,1){6.00}}
\put(85.00,18.00){\line(0,1){6.00}}
\put(100.00,18.00){\line(0,1){6.00}}

\put(70.00,24.00){\circle*{1.00}}
\put(75.00,28.00){\circle*{1.00}}
\put(65.00,28.00){\circle*{1.00}}

\put(80.00,28.00){\circle*{1.00}}
\put(85.00,24.00){\circle*{1.00}}
\put(90.00,28.00){\circle*{1.00}}

\put(95.00,28.00){\circle*{1.00}}
\put(100.00,24.00){\circle*{1.00}}
\put(105.00,28.00){\circle*{1.00}}

\put(70.00,18.00){\circle*{1.60}}
\put(85.00,18.00){\circle*{1.00}}
\put(100.00,18.00){\circle*{1.00}}

\qbezier(65,28)(65,24)(70,24)
\qbezier(70,24)(75,24)(75,28)

\qbezier(80,28)(80,24)(85,24)
\qbezier(85,24)(90,24)(90,28)
\qbezier(95,28)(95,24)(100,24)
\qbezier(100,24)(105,24)(105,28)


\put(51.00,-1.00){\scriptsize ${\mathcal G}_{1}\vartriangleright_{f} {\mathcal G}_{2}:$}
\put(69.00,-4.00){\scriptsize $e$}

\put(80.00,4.00){\circle*{1.00}}
\put(85.00,-6.00){\circle*{1.00}}
\put(90.00,4.00){\circle*{1.00}}

\put(95.00,4.00){\circle*{1.00}}
\put(100.00,-6.00){\circle*{1.00}}
\put(105.00,4.00){\circle*{1.00}}

\put(70.00,0.00){\circle*{1.60}}
\put(85.00,0.00){\circle*{1.00}}
\put(100.00,0.00){\circle*{1.00}}

\put(70.00,0.00){\line(1,0){30.00}}
\put(70.00,0.00){\line(0,1){6.00}}
\put(85.00,-6.00){\line(0,1){6.00}}
\put(100.00,-6.00){\line(0,1){6.00}}

\put(70.00,6.00){\circle*{1.00}}
\put(75.00,10.00){\circle*{1.00}}
\put(65.00,10.00){\circle*{1.00}}

\qbezier(65,10)(65,6)(70,6)
\qbezier(70,6)(75,6)(75,10)

\qbezier(80,4)(80,0)(85,0)
\qbezier(85,0)(90,0)(90,4)
\qbezier(95,4)(95,0)(100,0)
\qbezier(100,0)(105,0)(105,4)
\end{picture}

\caption{The comb product of rooted graphs
${\mathcal G}_{1}\vartriangleright {\mathcal G}_{2}=({\mathcal G}_{1},e_1)\vartriangleright ({\mathcal G}_{2},e_2)$ 
versus the essential component 
${\mathcal G}_{1}\vartriangleright_{f_2}{\mathcal G}_{2}$ of the 
c-comb product of the associated birooted graphs
in the case when $f_2\neq e_2$.
In the first product we only use $e_2$ as the glueing vertex 
for ${\mathcal G}_{2}$, whereas in the second product we use both $e_2$ and $f_2$.}
\end{figure}

\begin{Remark}
{\rm 
Note that instead of birooted graphs we could equivalently use two graphs and two roots in our constructions especially since 
the distributions with respect to the considered states depend only on the connected components of the birooted graphs. 
Observe also that using the simplified notations
\begin{eqnarray*}
{\mathcal G}_{1}\vartriangleright {\mathcal G}_{2}&=&({\mathcal G}_{1},f_1)\vartriangleright ({\mathcal  G}_{2},f_2)\\
{\mathcal G}_{1}\vartriangleright_{f_2} {\mathcal G}_{2}&=&({\mathcal G}_{1},e_1)\vartriangleright_{f_2} ({\mathcal  G}_{2},e_2),
\end{eqnarray*}
we can write Definition 4.2 as 
$$
{\mathcal G}=({\mathcal G}_{1}\vartriangleright_{f_2} {\mathcal G}_{2})
\oplus 
({\mathcal G}_{1}\vartriangleright {\mathcal G}_{2}).
$$ 
Of course, since the second component of this disjoint union is simply the comb product, all new information about the c-comb product is contained in the essential component. We use disconnected graphs consisting of two connected components in order to conform with the definition of the c-comb product of states of Hasebe.}
\end{Remark}

\begin{Theorem}
Let $({\mathcal  G}_{1},e_{1},f_1)$ and $({\mathcal  G}_{2},e_{2},f_2)$ be 
birooted graphs with adjacency matrices $a_1$ and $a_2$, respectively, and let 
$f_2\in V_2$. The adjacency matrix of ${\mathcal G}={\mathcal  G}_{1}\vartriangleright_{f_2} {\mathcal  G}_{2}$ has a decomposition
$$
a({\mathcal G}) =S_1+S_2,
$$
where the summands take the form
\begin{eqnarray*}
S_1&=&a_1\otimes P_{e_2}\otimes  P_{f_2}\\
S_2&=&P_{e_1}\otimes a_2 \otimes 1_{2}
+P_{e_1}^{\perp}\otimes 1_{2}\otimes a_{2},
\end{eqnarray*}
and act on the Hilbert space 
${\mathcal H}=\ell^2(V_1\vartriangleright_{f_2} V_2)\subset \ell^{2}(V_1)\otimes \ell^{2}(V_2)\otimes \ell^{2}(V_{2})$, with 
the root $e$ identified with $e_1\otimes e_2\otimes f_2$.
\end{Theorem}
{\it Proof.}
It follows from Definition 4.1 that the construction of the c-comb product of graphs can be carried out in two steps: 
\begin{enumerate}
\item[(a)]
take the orthogonal product of $({\mathcal  G}_{1},e_1)$ 
and $({\mathcal  G}_{2},f_2)$, obtained by attaching ${\mathcal G}_{2}$ at the vertex $f$ to all vertices but the root of ${\mathcal G}_{1}$, \item[(b)]
take the star product of this rooted graph with $({\mathcal G}_{2}, e_2)$ by glueing both graphs at their roots.
\end{enumerate}
In view of this decomposition, it is convenient 
to treat the set of vertices $V$ of ${\mathcal  G}_{1}\vartriangleright_{f_2} {\mathcal  G}_{2}$ as a subset of
$V_1\times V_2\times V_2$ of the 3d-space and the adjacency matrix as an operator living in 
the subspace of $\ell^{2}(V_1)\otimes \ell^{2}(V_2)\otimes \ell^{2}(V_2)$ of the form
\begin{eqnarray*}
{\mathcal H}&=&
({\mathcal H}_{1}\otimes f_2\otimes e_2)\oplus ({\mathcal H}_{1}^{\circ}\otimes \widetilde{{\mathcal H}}_{2}^{\circ}\otimes e_{2})
\oplus (e_1\otimes f_2\otimes {\mathcal H}_{2}^{\circ})\\
&\cong& ({\mathcal H}_{1}\otimes e_2\otimes f_2)\oplus 
({\mathcal H}_{1}^{\circ}\otimes e_{2}\otimes \widetilde{{\mathcal H}}_{2}^{\circ})
\oplus (e_1\otimes {\mathcal H}_{2}^{\circ}\otimes f_2)\\
&\cong& {\mathcal H}_{1}\oplus {\mathcal H}_{2}^{\circ} \oplus ({\mathcal H}_{1}^{\circ}\otimes 
\widetilde{{\mathcal H}}_{2}^{\circ})
\end{eqnarray*}
where ${\mathcal H}_{1}^{\circ}=\ell^{2}(V_1^{\circ})$, ${\mathcal H}_{2}^{\circ}=\ell^{2}(V_2^{\circ})$
and $\widetilde{{\mathcal H}}_{2}^{\circ}=\ell^{2}(V_{2}\setminus \{f\})$.
For convenience, in the first step we 
applied the flip $\Sigma_{2,3}(a\otimes b \otimes c)=a\otimes c \otimes b$.
Using the above decomposition of the c-comb product, we can construct the adjacency matrix of ${\mathcal  G}_{1}\vartriangleright_{f_2} {\mathcal  G}_{2}$ 
from those of the orthogonal and star products of rooted graphs:
\begin{eqnarray*}
a({\mathcal  G}_{1}\vartriangleright_{f_2} {\mathcal  G}_{2})&=&
a(({\mathcal  G}_{1},e_1)\vdash ({\mathcal  G}_{2},f_2))\otimes P_{e_2}+P_{(e_1,f_2)}\otimes a_2\\
&=&
(a_1\otimes P_{f_2}+P_{e_1}^{\perp}\otimes a_2)\otimes P_{e_2}+P_{e_1}\otimes 
P_{f_2}\otimes a_2\\
&\cong&
a_1\otimes P_{e_2}\otimes P_{f_2}+P_{e_1}\otimes a_2\otimes P_{f_2}+P_{e_1}^{\perp}\otimes P_{e_2}\otimes a_2,
\end{eqnarray*}
where the last step corresponds to the flip $\Sigma_{2,3}$, which gives a more convenient form, 
at least from the point of view of Proposition 3.1 and earlier works on independence [11,12]. 
After applying the flip, 
the root $e$ is identified with $e_1\otimes e_2 \otimes f_2$. Now, it remains to observe that
\begin{eqnarray*}
P_{e_1}\otimes a_2\otimes P_{f_2}&\equiv& P_{e_1}\otimes a_2\otimes 1_{2}\\
P_{e_1}^{\perp}\otimes P_{e_2}\otimes a_2 &\equiv& P_{e_1}^{\perp}\otimes 1_{2}\otimes a_2
\end{eqnarray*}
on ${\mathcal H}$ (they act on ${\mathcal H}$ identically). This completes the proof.
\hfill $\blacksquare$\\

\begin{Corollary}
The adjacency matrix of the c-comb product 
${\mathcal G}=({\mathcal G}_{1},e_1,f_1)\vartriangleright 
({\mathcal G}_{2},e_{e},f_{2})$ is of the form
$$
a({\mathcal  G})=
a({\mathcal G}_{1}\vartriangleright_{f_2} {\mathcal G}_{2})\oplus
a({\mathcal G}_{1}\vartriangleright {\mathcal G}_{2}) 
$$
in the Hilbert space ${\mathcal H}\oplus {\mathcal H}'$, equipped with the canonical inner product, where
$$
a({\mathcal  G}_{1}\vartriangleright {\mathcal  G}_{2}) =
a_1\otimes  P_{f_2}+1_1\otimes a_{2}
$$
is the adjacency matrix of ${\mathcal  G}_{1}\vartriangleright{\mathcal  G}_{2}$ 
in ${\mathcal H}'=\ell^{2}(V_1)\otimes \ell^{2}(V_{2})$.
\end{Corollary}
{\it Proof.}
The form of $a({\mathcal G})$ follows from Definition 4.2 and 
$a({\mathcal  G}_{1}\vartriangleright {\mathcal  G}_{2})$ was given in 
Remark 4.1(3), which proves our claim.
\hfill $\blacksquare$\\

\begin{Corollary}
With the notations of Corollary 4.1, $a({\mathcal G})=S_1+S_2$, where 
\begin{eqnarray*}
S_{1}&=&a_{1}\otimes P_{e_2}\otimes P_{f_2}+a_1\otimes P_{f_2}\\
S_{2}&=&P_{e_1}\otimes a_{2}\otimes 1_{2}
+P_{e_1}^{\perp}\otimes 1_{2}\otimes a_{2}+1\otimes a_{2}
\end{eqnarray*}
and the pair $(S_1,S_2)$ is c-monotone independent with respect to 
$(\varphi, \psi)$, where 
$$
\varphi=\widetilde{\varphi}_1\otimes \widetilde{\varphi}_2\otimes \widetilde{\psi}_2,\;\;\;
\psi=\widetilde{\psi}_1\otimes \widetilde{\psi}_2
$$
and $\varphi_1$, $\varphi_2$, $\psi_2$ are states associated with $e_1\in V_1$ and $e_2,f_2\in V_2$, 
respectively.
\end{Corollary}
{\it Proof.} 
The equation $a({\mathcal G})=S_1+S_2$
follows from Theorem 4.1 and Corollary 4.1,
Moreover, the joint distribution of $(S_1,S_2)$ agrees with the joint
distribution of $(S_1',S_2')$ in the state $\varphi$, where
\begin{eqnarray*}
S_{1}'&=&a_1\otimes P_{f_2}\otimes P_{e_2}\\
S_{2}'&=&P_{e_1}\otimes a_2\otimes 1_{2}
+P_{e_1}^{\perp}\otimes 1_{2}\otimes a_{2}
\end{eqnarray*}
since ${\mathcal H}\perp {\mathcal H}'$, where ${\mathcal H}, {\mathcal H}'$
are as in Theorem 4.1 and Corollary 4.1. 
However, by Proposition 3.1, the joint distribution of $S_1',S_2'$ in the 
state $\varphi$ agrees with that of a pair of c-monotone random variables
w.r.t. the state $\varphi$ (the projections $P_{e_1},P_{e_2}$ and $P_{f_2}$ are 
images of the GNS representations of generic projections $P,Q$ of Proposition 3.1 when
we take noncommutative probability spaces defined by extended states $\widetilde{\varphi}_{1}, \widetilde{\varphi}_{2}, 
\widetilde{\psi}_{1},\widetilde{\psi}_{2}$, described in Section 3).
Moreover, the joint distribution of $S_1,S_2$ in the state $\psi$ agrees with that of $S_{1}''= a_1\otimes P_{f_{2}}$ and 
$S_{2}''=1_{1}\otimes a_2$ since $S_1',S_2'$ act trivially on ${\mathcal H}'$.
Therefore, the pair $(S_1,S_2)$ is c-monotone independent w.r.t. $(\varphi, \psi)$.
This completes the proof.
\hfill $\blacksquare$\\

\begin{Proposition}
Let $\mu_1, \mu_2,\nu_2$ be the spectral distributions of $({\mathcal G}_{1},e_1)$, 
$({\mathcal G}_{2},e_2)$ and $({\mathcal G}_{2}, f_2)$, respectively.
Then the spectral distribution of $({\mathcal  G}_{1}\vartriangleright_{f_2} {\mathcal  G}_{2},e)$ is given by 
the c-monotone additive convolution $\mu=\mu_1 \vartriangleright_{\nu_2} \mu_2$, where 
the reciprocal Cauchy transform of this distribution is given by
$$
F_{\mu}(z)=F_{\mu_{1}}(F_{\nu_{2}}(z))+F_{\mu_2}(z)-F_{\nu_2}(z)
$$
where $F_{\mu}(z)=1/G_{\mu}(z)$ and $G_{\mu}(z)=\sum_{n=0}^{\infty}M_{n}z^{-n-1}$, 
where $M_{n}$ is the $n$-th moment of $\mu$ (equal to the number of 
walks of lenght $n$ from the root to the root).
\end{Proposition}
{\it Proof.}
By Theorem 4.1, the spectral distribution of the c-comb product of graph is given by 
the c-monotone additive convolution $\mu_1\vartriangleright_{\nu_2}\mu_2$. 
Thus, we can use the formula for the reciprocal Cauchy 
transform of $\mu:=\mu_{1}\vartriangleright_{\nu_2}\mu_2$ derived by Hasebe [8] which agrees 
with the above formula. Let us present another short proof, based on 
Definition 4.1:
\begin{eqnarray*}
F_{\mu}(z)&=&F_{\mu_2}(z)+F_{\mu_1\vdash\nu_2}(z)-z\\
&=& F_{\mu_2}(z)+ (F_{\mu_{1}}(F_{\nu_2}(z))-F_{\nu_2}(z)+z)-z\\
&=& F_{\mu_{1}}(F_{\nu_2}(z))+F_{\mu_2}(z) -F_{\nu_2}(z)
\end{eqnarray*}
and that completes the proof.
\hfill $\blacksquare$\\

\begin{Corollary}
With the notations of Proposition 4.2, the spectral distribution of the c-comb product of graphs $({\mathcal  G},e,f)$ is the pair $(\mu, \nu)$, 
where $\mu= \mu_1 \vartriangleright_{\nu_2}\mu_2$ and 
$\nu=\nu_1\vartriangleright \nu_2$.
\end{Corollary}
{\it Proof.}
The spectral distribution of a birooted graph is a pair of spectral distributions associated with both roots. Since ${\mathcal G}$ is a disjoint sum of two connected graphs, it suffices to determine the spectral distributions of these components. By Proposition 4.2, the essential component has the spectral distribution given by $\mu$ and the second component has the distribution given by $\nu$. This ends the proof.\hfill $\blacksquare$

\section{C-comb loop product of graphs}

We have already seen that the spectral distributions of the c-comb product of rooted graphs is given by the c-monotone additive 
convolution of probability measures. A similar correspondence holds for the additive convolutions associated with other types of noncommutative independence (free, boolean, orthogonal, s-free) and the associated products of rooted graphs (free, star, orthogonal, s-free). 

We have shown in [14] that a similar correspondence holds for the multiplicative convolutions associated with these types
of noncommutative independence except that one has to introduce {\it loop products} of rooted graphs. Roughly speaking, 
each of the loop products is obtained from the original product by adding loops to appropriate vertices. 
This operation corresponds to the observation made by Bercovici [4] that that in order to introduce a multiplicative convolution associated with monotone independence, one needs to modify the usual monotone independent variables by adding a certain projecton (for the boolean multiplicative convolution, see Franz [7]). Then, a nice combinatorial formula can be derived for the moments of the product of the modified 
adjacency matrices.

\begin{Remark}
{\rm 
Let us recall basic facts on the loop products of rooted graphs introduced and studied in [14].
\begin{enumerate}
\item
In order to construct a loop product of rooted 
graphs, $({\mathcal G}_{1},e_1)$ and $({\mathcal G}_{2},e_2)$, associated with 
some noncommutative independence (boolean, free, monotone, s-free, orthogonal), we first construct the usual 
product of graphs associated with this independence  (star, free, comb, s-free, orthogonal). 
Suppose the edges of the product graph are {\it naturally colored}, by which we understand that all edges of all copies of ${\mathcal G}_{1}$ are colored by 1 and all edges of all copies of ${\mathcal G}_{2}$ are colored by 2.
By an {\it alternating walk} on the product graph we understand a walk in which consecutive edges have different colors, 
where all edges of copies of ${\mathcal G}_{i}$ are colored by 
$i\in \{1,2\}$. We assume that the rooted graphs are locally finite.
It is convenient (but not necessary) to consider rooted graphs which have loops at their roots 
to maximize the number of alternating root-to-root alternating walks.
\item
The {\it comb loop product} of rooted graphs $({\mathcal G}_{1},e_1)$ and $({\mathcal G}_{2},e_2)$
is the rooted graph $({\mathcal G},e)=({\mathcal  G}_{1}\vartriangleright_{\ell} {\mathcal  G}_{2},e)$ 
obtained from the comb product 
$({\mathcal G}_{1}\vartriangleright {\mathcal G}_{2},e)$
by attaching a loop of color $1$ to each vertex but the root of each copy of ${\mathcal G}_{2}$. 
Its adjacency matrix has a decomposition
$$
a({\mathcal G})=R_1+R_2,
$$
where the one-color adjacency matrices are of the form
$$
R_1=a_1\otimes P_{e_2}+ 1_1\otimes P_{e_2}^{\perp} \;\;\;{\rm and}\;\;\; R_2=1_1\otimes a_2,
$$
for $a_{i},1_{i}\in B(\ell^{2}(V_{i}))$, with the identity operators
$1_{i}$, and $i=1,2$. The pair $(R_1-1,R_2-1)$, where $1=1_1\otimes 1_2$, is monotone independent w.r.t.
the vector state $\varphi_{e}$ associated with $\delta(e)$. Thus the moments of the product $R_{2}R_{1}$ agree with the moments of the multiplicative monotone convolution $\mu_1\circlearrowright \mu_2$.  
Note that $R_1$ is obtained from the usual comb product adjacency matrix $S_{1}$ by adding 
a projection $L_{1}=1_1\otimes P_{e_2}^{\perp}$, whereas $R_{2}=S_{2}$. 
\item
A convenient tool used in the study of multiplicative convolutions is the so-called $\psi$-transform.
Let ${\mathcal  M}_{{\mathbb R}_{+}}$ denote the set
of probability measures on ${\mathbb R}_{+}=[0, \infty)$.
If $\mu\in{\mathcal M}_{{\mathbb R}_{+}}$, we can define
$$
\psi_{\mu}(z)=\int_{{\mathbb R}_{+}}\frac{zt}{1-zt}d\mu(t), \;\;\; z\in {\mathbb C}\setminus 
{\mathbb R}_{+},
$$
which, in the case when $\mu$ has finite moments of all orders, becomes the 
moment generating function $\psi_{\mu}(z)=\sum_{n=1}^{\infty}\mu(X^n)z^{n}$, 
where $\mu(X^n)$ are the moments
of the unique functional $\mu:{\mathbb C}[X]\rightarrow {\mathbb C}$ defined by $\mu$.
\item
A central role is played by the transform related to $\psi_{\mu}(z)$, namely 
$$
\eta_{\mu}(z)=\frac{\psi_{\mu}(z)}{1+\psi_{\mu}(z)}=\sum_{n=1}^{\infty}N_{\mu}(n)z^{n}
$$
where $z\in {\mathbb C}\setminus {\mathbb R}_{+}$. More 
generally, if $Z$ is any random variable, then we denote by
$\psi_Z$ the moment generating function (without the constant term) as a formal power series, and 
by $\eta_{Z}$ the formal power series given by the above formula. If $Z$ is an adjacency matrix, the coefficients $N_{Z}(n)$ of $\eta_{Z}$ can be interpreted as the numbers of root-to-root first return walks, or {\it f-walks}.
\item
In order to formulate our multiplicative results, we take now $\eta_{Z}$ 
for $Z=R_{2}R_{1}$ in the state $\varphi_{e}$ associated with the root $e$, namely
$$
\eta_{Z}(z)=\sum_{n=1}^{\infty}N_{Z}(n)z^{n}.
$$
Since $Z$ is a product (of one-color adjacency matrices), the corresponding walks on the product graph 
must have alternating colors and they must be of even lenght. Moreover, due to 
the way the loops are added to all product graphs, the `first return moments' for
$Z$ are associated with {\it rooted alternating d-walks} (by a {\it d-walk} we understand 
a `double return walk originating with color 1') 
counted on different products.
Enumeration results for rooted alternating d-walks were obtained in [14].
\item
In particular, if ${\mathcal G}_{1}{\vartriangleright}_{\ell}\,{\mathcal G}_{2}$ is naturally colored
and $A({\mathcal G}_{1}{\vartriangleright}_{\ell}\,{\mathcal G}_{2})=R_1+R_2$ is the
decomposition of its adjacency matrix induced by the coloring, then 
$$
N_{Z}(n)=N_{\mu_1\circlearrowright\mu_2}(n)=|D_{2n}(e)|
$$
where $Z=R_2R_1$ and $D_{2n}(e)$ denotes the set of rooted alternating d-walks 
on ${\mathcal G}_{1}{\vartriangleright}_{\ell}\,{\mathcal G}_{2}$ 
of length $2n$, where $n\in {\mathbb N}$. 
Here, $\mu_1\circlearrowright \mu_2$ is the monotone multiplicative convolution of $\mu_1$ and $\mu_2$,
defined by the formula
$$
\eta_{\mu_1\circlearrowright \mu_2}(z)=\eta_{\mu_1}(\eta_{\mu_2}(z)),
$$
for any $\mu_1,\mu_2\in {\mathcal M}_{{\mathbb R}^{+}}$, using the original notation of Bercovici [4] for the monotone multiplicative convolution [4], different than the new notation of Hasebe [9].
Analogous results were established in [14] for star, orthogonal, free and s-free products of graphs
and the associated convolutions.
\end{enumerate}
}
\end{Remark}

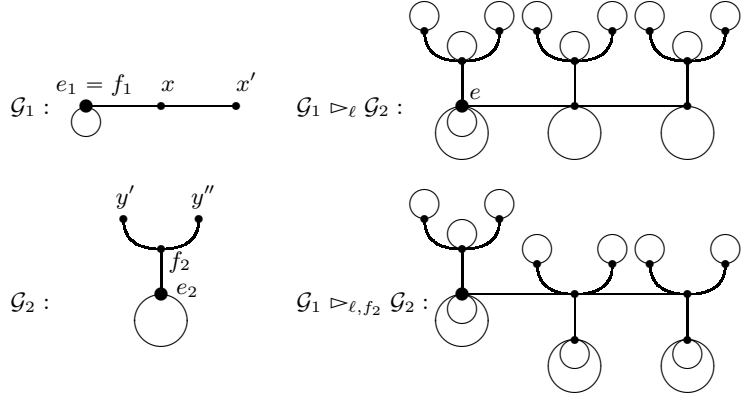
\begin{figure}
\unitlength=1mm
\special{em.linewidth 2pt}
\linethickness{0.5pt}
\begin{picture}(80.00,60.00)(20.00,-23.00)
\put(20.00,15.00){\line(1,0){20.00}}

\put(20.00,15.00){\circle*{1.60}}
\put(20.00,13.00){\circle{4.00}}
\put(30.00,15.00){\circle*{1.00}}
\put(40.00,15.00){\circle*{1.00}}

\put(16.00,17.00){\scriptsize $e_1=f_1$}
\put(71.00,16.00){\scriptsize $e$}
\put(30.00,17.00){\scriptsize $x$}
\put(40.00,17.00){\scriptsize $x'$}
\put(10.00,14.00){\scriptsize ${\mathcal G}_{1}:$}

\put(30.00,-10.00){\line(0,1){6.00}}

\put(30.00,-10.00){\circle*{1.60}}
\put(30.00,-13.50){\circle{7.00}}
\put(30.00,-4.00){\circle*{1.00}}
\put(25.00,0.00){\circle*{1.00}}
\put(35.00,0.00){\circle*{1.00}}

\put(32.00,-10.00){\scriptsize $e_{2}$}
\put(31.00,-6.50){\scriptsize $f_2$}
\put(34.00,2.00){\scriptsize $y''$}
\put(24.00,2.00){\scriptsize $y'$}
\put(10.00,-12.00){\scriptsize ${\mathcal G}_{2}:$}

\qbezier(25,0)(25,-4)(30,-4)
\qbezier(30,-4)(35,-4)(35,0)


\put(65.00,27.00){\circle{4.00}}
\put(70.00,23.00){\circle{4.00}}
\put(75.00,27.00){\circle{4.00}}

\put(80.00,27.00){\circle{4.00}}
\put(85.00,23.00){\circle{4.00}}
\put(90.00,27.00){\circle{4.00}}

\put(95.00,27.00){\circle{4.00}}
\put(100.00,23.00){\circle{4.00}}
\put(105.00,27.00){\circle{4.00}}

\put(48.00,14.00){\scriptsize ${\mathcal G}_{1}\vartriangleright_{\ell} {\mathcal G}_{2}:$}

\put(70.00,15.00){\line(1,0){30.00}}
\put(70.00,15.00){\line(0,1){6.00}}
\put(85.00,15.00){\line(0,1){6.00}}
\put(100.00,15.00){\line(0,1){6.00}}

\put(70.00,13.00){\circle{4.00}}
\put(70.00,11.50){\circle{7.00}}
\put(85.00,11.50){\circle{7.00}}
\put(100.00,11.50){\circle{7.00}}

\put(70.00,21.00){\circle*{1.00}}
\put(75.00,25.00){\circle*{1.00}}
\put(65.00,25.00){\circle*{1.00}}

\qbezier(65,25)(65,21)(70,21)
\qbezier(70,21)(75,21)(75,25)

\qbezier(80,25)(80,21)(85,21)
\qbezier(85,21)(90,21)(90,25)
\qbezier(95,25)(95,21)(100,21)
\qbezier(100,21)(105,21)(105,25)

\put(70.00,15.00){\circle*{1.60}}
\put(85.00,15.00){\circle*{1.00}}
\put(100.00,15.00){\circle*{1.00}}

\put(85.00,21.00){\circle*{1.00}}
\put(100.00,21.00){\circle*{1.00}}
\put(90.00,25.00){\circle*{1.00}}
\put(80.00,25.00){\circle*{1.00}}
\put(105.00,25.00){\circle*{1.00}}
\put(95.00,25.00){\circle*{1.00}}


\put(65.00,02.00){\circle{4.00}}
\put(70.00,-02.00){\circle{4.00}}
\put(75.00,02.00){\circle{4.00}}

\put(80.00,-04.00){\circle{4.00}}
\put(85.00,-18.00){\circle{4.00}}
\put(90.00,-04.00){\circle{4.00}}

\put(95.00,-04.00){\circle{4.00}}
\put(100.00,-18.00){\circle{4.00}}
\put(105.00,-04.00){\circle{4.00}}

\put(48.00,-12.00){\scriptsize ${\mathcal G}_{1}\vartriangleright_{\ell, f_2} {\mathcal G}_{2}:$}

\put(70.00,-10.00){\line(1,0){30.00}}
\put(70.00,-10.00){\line(0,1){6.00}}
\put(85.00,-16.00){\line(0,1){6.00}}
\put(100.00,-16.00){\line(0,1){6.00}}

\put(70.00,-12.00){\circle{4.00}}
\put(70.00,-13.50){\circle{7.00}}
\put(85.00,-19.50){\circle{7.00}}
\put(100.00,-19.50){\circle{7.00}}

\put(70.00,-04.00){\circle*{1.00}}
\put(75.00,0.00){\circle*{1.00}}
\put(65.00,0.00){\circle*{1.00}}

\qbezier(65,0)(65,-4)(70,-4)
\qbezier(70,-4)(75,-4)(75,0)

\qbezier(80,-6)(80,-10)(85,-10)
\qbezier(85,-10)(90,-10)(90,-6)

\qbezier(95,-6)(95,-10)(100,-10)
\qbezier(100,-10)(105,-10)(105,-6)

\put(70.00,-10.00){\circle*{1.60}}
\put(85.00,-16.00){\circle*{1.00}}
\put(100.00,-16.00){\circle*{1.00}}

\put(85.00,-10.00){\circle*{1.00}}
\put(100.00,-10.00){\circle*{1.00}}
\put(90.00,-6.00){\circle*{1.00}}
\put(80.00,-6.00){\circle*{1.00}}
\put(105.00,-6.00){\circle*{1.00}}
\put(95.00,-6.00){\circle*{1.00}}

\end{picture}
\caption{The comb loop product of rooted graphs 
${\mathcal G}_{1}\vartriangleright_{\ell}{\mathcal G}_{2}=({\mathcal G}_{1},e_1)\vartriangleright_{\ell}({\mathcal G}_{2},e_2)$ 
versus the essential component ${\mathcal G}_{1}\vartriangleright_{\ell, f_2}\,{\mathcal G}_{2}$ of the 
c-comb loop product of birooted graphs 
$({\mathcal G}_{1},e_1,f_1)\vartriangleright_{\ell} ({\mathcal G}_{2},e_2,f_2)$.}
\end{figure}

In the case of c-monotone independence, the right category of graphs is that of birooted graphs. Again, 
we assume that they are uniformly locally finite, which allows us to view their adjacency matrices as bounded operators on 
$\ell^{2}(V)$. We would like to construct a loop product of birooted graphs associated with c-monotone independence. 
As in the case of loop products of rooted graphs studied in [14], we will assume that all edges of all copies of ${\mathcal G}_{1}$ and ${\mathcal G}_{2}$ are colored by $1$ and $2$, respectively. The idea of the loop product of rooted graphs construction for ${\mathcal I}$-independence is based on adding loops to the usual product of rooted graphs for ${\mathcal I}$-independence
in such a way that the spectral distribution of the product of their adjacency matrices is given by the 
multiplicative convolution associated with ${\mathcal I}$.
In the case of c-monotone independence, we have a pair of measures, and thus this convolution is a pair 
of convolutions.

\begin{Definition}
{\rm  Let $({\mathcal  G}_{1},e_1,f_1)$ and $({\mathcal  G}_{2},e_2,f_2)$ be birooted graphs. By their {\it c-comb loop product} denoted 
$({\mathcal G}_{1},e_1,f_1)\vartriangleright_{\ell}
({\mathcal G}_{2},e_2,f_2)$ we understand the disjoint union of rooted graphs 
$$
({\mathcal G},e,f)=\left(({\mathcal G}_{1},e_1)\vartriangleright_{\ell,f_2}({\mathcal G}_{2},e_2)\right)\oplus 
\left(({\mathcal G}_{1},f_1)\vartriangleright_{\ell} ({\mathcal G}_2,f_2)\right),
$$
where $({\mathcal  G}_{1},e_1)\vartriangleright_{\ell,f_2}({\mathcal  G}_{2},e_2)$, 
called the {\it essential component} of $({\mathcal G},e,f)$, is the rooted graph obtained from 
$({\mathcal G}_{1},e_1)\vartriangleright_{f_2}({\mathcal G}_{2},e_2)$ by adding a loop of color 2 to each vertex 
but $e_2$ of the copy of ${\mathcal G}_{2}$ attached to $e_1$ and to each vertex but $f_2$ 
of all the remaining copies of ${\mathcal G}_{2}$. The roots $e,f$ are identified with $(e_1,f_2,e_2)$ and $(f_1,f_2)$, 
respectively (cf. Definition 4.2). 
We will also use the simplified notation 
${\mathcal G}_{1}\vartriangleright_{\ell, f_2}{\mathcal G}_{2}$ for the essential component. 
An example is given in Fig. 2.
}
\end{Definition}

\begin{Proposition}
Let $({\mathcal  G}_{1},e_{1},f_1)$ and $({\mathcal  G}_{2},e_{2},f_2)$ be 
birooted graphs with adjacency matrices $a_1$ and $a_2$, respectively. Then the adjacency matrix of 
${\mathcal  G}_{1}\vartriangleright_{\ell, f_2} {\mathcal  G}_{2}$ has the decomposition 
$$
a({\mathcal  G}_{1}\vartriangleright_{\ell, f_2} {\mathcal  G}_{2})=R_1+R_2,
$$
where the summands are given by formulas
\begin{eqnarray*}
R_1-1&=&a_1^{\curlyvee}\otimes P_{e_2}\otimes  P_{f_2}\\
R_2-1&=&P_{e_1}\otimes a_2^{\curlyvee} \otimes 1_{2}
+P_{e_1}^{\perp}\otimes 1_{2}\otimes a_{2}^{\curlyvee}
\end{eqnarray*}
and act on vectors from the Hilbert space ${\mathcal H}=\ell^{2}(V_{1}\vartriangleright _{f_2}V_{2})$, where
$a_i^{\curlyvee}=a_i-1_i$ for $i=1,2$. In this realization, 
the root $e$ of ${\mathcal G}$ is identified with $e_1\otimes e_2\otimes f_2$. 
\end{Proposition}
{\it Proof.}
Since $1=1_1\otimes 1_2\otimes 1_2$, we can use the decompositions 
$$
1_2=P_{f_2}+P_{f_2}^{\perp}=P_{e_{2}}+P_{e_{2}}^{\perp}
$$
to write $R_1$ and $R_2$ in a different form:
$$
R_1=S_1+1_1\otimes 1_2\otimes P_{f_2}^{\perp}+1_1\otimes P_{e_2}^{\perp}\otimes P_{f_2}\;\;\;{\rm and}\;\;\;
R_2=S_2
$$
where $S_1,S_2$ are as in Theorem 4.1. Under the identification of the roots 
of ${\mathcal  G}_{1}\vartriangleright_{f_2} {\mathcal  G}_{2}$ (as in Definition 4.1) with basis 
vectors of the Hilbert space $\ell^{2}(V_1\vartriangleright_{f_2}V_2)$, the following 
correspondences between projections in the expression for $R_1$ and loops take place:
\begin{enumerate}
\item
the projection $1_1\otimes 1_2\otimes P_{f_2}^{\perp}\equiv 1_{1}\otimes P_{e}\otimes P_{f_2}^{\perp}$ corresponds to 
loops of color 1 added to all vertices but $f_2$ in all copies of ${\mathcal G}_{2}$ which are not attached to
the root $e_1$ of ${\mathcal G}_{1}$, where $P_1\equiv P_2$ means that 
projections $P_1$ and $P_2$ act identically on the underlying Hilbert space ${\mathcal H}$, 
which can be seen by using the direct sum decomposition of ${\mathcal H}$ given in the proof of Theorem 4.1,
\item
the projection 
$1_1\otimes P_{e_2}^{\perp}\otimes P_{f_2}$ corresponds to the loops of color 1 added to all vertices
but $e_2$ in the copy of ${\mathcal G}_{2}$ attached to the root $e_1$ of ${\mathcal G}_{1}$.
\end{enumerate}
In order to visualize this correspondence geometrically, one should identify the vertices of $V_{1}\vdash V_{2}$ with vectors
of the 2d plane formed by the first and third legs of the tensor product and the vertices of the second $V_2$ 
in the set $(V_{1}\vdash V_2)\star V_2$ with the vectors that appear at the second leg.
This proves that $a({\mathcal G})=R_1+R_2$, which completes the proof.
\hfill $\blacksquare$\\

\begin{Lemma}
Let $({\mathcal  G}_{1},e_{1},f_1)$ and $({\mathcal  G}_{2},e_{2},f_2)$ be 
birooted graphs with adjacency matrices $a_1$ and $a_2$, respectively, and 
let $({\mathcal G},e,f)=({\mathcal G}_{1},e_1,f_1)\vartriangleright_{\ell}
({\mathcal G}_{2},e_{2},f_{2})$. Then
\begin{enumerate}
\item
the adjacency matrix of $({\mathcal G},e,f)$ has the decomposition $a({\mathcal G})=R_1+R_2$, where 
\begin{eqnarray*}
R_{1}-1&=&(a_{1}^{\curlyvee}\otimes P_{e_2}\otimes P_{f_2})\oplus (a_1^{\curlyvee}\otimes P_{f_2})\\
R_{2}-1&=&(P_{e_1}\otimes a_{2}^{\curlyvee}\otimes 1_{2}
+P_{e_1}^{\perp}\otimes 1_{2}\otimes a_{2}^{\curlyvee})\oplus (1_1\otimes a_{2}^{\curlyvee}),
\end{eqnarray*}
acting on vectors from the Hilbert space ${\mathcal H}\oplus {\mathcal H}'$, 
\item
the pair $(R_1-1,R_2-1)$ is c-monotone independent with respect to 
$(\varphi, \psi)$, where 
$$
\varphi=\widetilde{\varphi}_1\otimes \widetilde{\varphi}_2\otimes \widetilde{\psi}_2,\;\;\;
\psi=\widetilde{\psi}_1\otimes \widetilde{\psi}_2
$$
and $\widetilde{\varphi}_i, \widetilde{\psi}_i$ are states associated with $e_i,f_i\in V_i$, where $i=1,2$, 
respectively.
\end{enumerate}
\end{Lemma}
{\it Proof.}
By Definition 5.1, the adjacency matrix of ${\mathcal G}$ has the form
$$
a({\mathcal  G})=
a({\mathcal G}_{1}\vartriangleright_{\ell, f_2} {\mathcal G}_{2})\oplus
a({\mathcal G}_{1}\vartriangleright_{\ell} {\mathcal G}_{2}) 
$$
and acts in the Hilbert space ${\mathcal H}\oplus {\mathcal H}'$, equipped with the canonical inner product, where
$a({\mathcal G}_{1}\vartriangleright_{\ell, f_2} {\mathcal G}_{2})$ is given by Proposition 5.1 and
$$
a({\mathcal  G}_{1}\vartriangleright_{\ell} {\mathcal  G}_{2}) =
a_1\otimes  P_{f_2}+1_{1}\otimes P_{f_2}^{\perp}+1_1\otimes a_{2}
$$
is the adjacency matrix of ${\mathcal  G}_{1}\vartriangleright_{\ell}{\mathcal  G}_{2}$ 
in ${\mathcal H}'=\ell^{2}(V_1)\otimes \ell^{2}(V_{2})$ by Remark 5.1(2).
Now, the unit $1$ is the identity on ${\mathcal H}\oplus {\mathcal H}'$, hence
$1=(1_1\otimes 1_2\otimes 1_2)\oplus (1_2\otimes 1_2)$. Therefore, if we put together the results of 
Proposition 5.1 and Remark 5.1(2), we obtain the formulas for $R_1$ and $R_2$. The fact that 
the pair $(R_1-1, R_2-1)$ is c-monotone independent follows then from Proposition 3.1. This completes the proof.
\hfill $\blacksquare$\\

We would like to establish a connection between walks on the c-comb product of graphs 
and the moments of the {\it c-monotone multiplicative convolution} of distributions on ${\mathbb C}[x]$, denoted 
$\Sigma$ (in particular, measures from ${\mathcal M}_{{\mathbb R}^{+}}$), defined by Hasebe in [9].
This convolution is a binary operation on the set $\Sigma \times \Sigma$, so the definition
involves pairs of distributions.
Nevertheless, we expect that a Multiplication Theorem similar to that proved in [13] for other
products of rooted graphs will also hold. 

For given $\mu_1,\nu_1,\mu_2,\nu_2\in \Sigma$, where $\nu_2$ is not concentrated at zero, let 
$$
(\mu_1,\nu_1)\circlearrowright (\mu_2, \nu_2)=(\mu_1\circlearrowright_{\nu_2}\mu_2, \nu_1\circlearrowright \nu_2),
$$ 
where 
\begin{eqnarray*}
\eta_{\mu_1\circlearrowright_{\nu_2}\mu_2}(z)
&=&
\frac{\eta_{\mu_2}(z)}{\eta_{\nu_{2}}(z)}\eta_{\mu_{1}}(\eta_{\nu_{2}}(z)).
\end{eqnarray*}
This convolution reminds the orthogonal multiplicative convolution $\mu_1\angle \nu_2$ [14], 
obtained when $\eta_{\mu_2}(z)=z$, which means that $\mu_2=\delta_{1}$.

\begin{Theorem}
Let $({\mathcal G},e,f)=({\mathcal G}_{1},e_1,f_1)\vartriangleright_{\ell}
({\mathcal G}_{2},e_{2},f_{2})$ be naturally colored
and let $A({\mathcal G})=R_1+R_2$ be the decomposition induced by the coloring of edges. Then
\begin{eqnarray*}
N_{Z,e}(n)&=&N_{\mu_1\circlearrowright_{\nu_2}\mu_2}(n)=|D_{2n}(e)|,\\
N_{Z,f}(n)&=&N_{\nu_1\circlearrowright \nu_2}(n)\;\;\;=|D_{2n}(f)|,
\end{eqnarray*}
where $Z=R_2R_1$ and $D_{2n}(e), D_{2n}(f)$ denote the sets of alternating 
d-walks of lenght $2n$ from e to e and from f to f, respectively, where $n\in {\mathbb N}$.
\end{Theorem}
{\it Proof.}
The second equation follows from that for the comb loop product of graphs [14]. Therefore, 
it suffices to show the first equation.
Using the formula for 
$\eta_{\mu_1\circlearrowright_{\nu_2}\mu_2}$ given above and those for $\eta$-series 
for the orthogonal multiplicative convolution $\mu_1\angle \mu_2$ and the boolean multiplicative 
convolution $\mu_1\boxasterisk\mu_2$, namely
$$
\eta_{\mu_1 \angle \mu_2}(z)=\frac{z\eta_{\mu_1}(\eta_{\mu_2}(z))}{\eta_{\mu_2}(z)}\;\;\;{\rm and}\;\;\;
\eta_{\mu_1\boxasterisk\,\mu_2}(z)=\frac{\eta_{\mu_1}(z)\eta_{\mu_2}(z)}{z},
$$ 
we obtain the decomposition
$$
\mu_1 \circlearrowright_{\nu_{2}}\mu_2=(\mu_1\angle \nu_2)\boxasterisk \mu_2.
$$ 
The coefficients of the $\eta$-series for the orthogonal loop product of rooted graphs are of the form
$$
N_{\mu_1 \angle \mu_2}(n)=
\sum_{r=2}^{n}
N_{\mu_1}(r)\sum_{k_{1}+k_{2}+\ldots +k_{r-1}=n-1}
N_{\mu_2}(k_{1})N_{\mu_2}(k_{2})\ldots
N_{\mu_2}(k_{r-1}),
$$
for $n\geq 2$, where it is assumed that the summation runs over positive integers 
$k_{1}, \ldots , k_{r-1}$ which add up to $n-1$ (one can formally include the term for $r=1$, as we did in [14], but in that case the second sum gives zero contribution). If $n=1$, we get $N_{\mu_1\angle \mu_2}(n)=N_{\mu_1}(1)$.
In turn, the coefficients of the $\eta$-series for the boolean loop product of graphs are of the form
$$
N_{\mu_1 \boxasterisk \,\mu_2}(n)=
\sum_{j+k=n+1}N_{\mu_1}(j)N_{\mu_2}(k)
$$
where the summation runs over positive integers which add up to $n+1$.
This gives 
\begin{eqnarray*}
N_{\mu_1 \circlearrowright_{\nu_{2}}\mu_2}(n)&=&
\sum_{j+k=n+1}
\left(\sum_{r=2}^{j}N_{\mu_1}(r	)
\sum _{k_{1}+k_{2}+\ldots +k_{r-1}=j-1}
N_{\nu_1}(k_1) \cdots
N_{\nu_2}(k_{r-1})\right)N_{\mu_2}(k)\\
&=&
\sum_{r=2}^{n}N_{\mu_1}(r)\sum _{k_{1}+k_{2}+\ldots +k_{r}=n}
N_{\nu_2}(k_{1})\ldots
N_{\nu_2}(k_{r-1})N_{\mu_2}(k_r)
\end{eqnarray*}
for $n\geq 2$, whereas $N_{\mu_1 \circlearrowright_{\nu_{2}}\mu_2}(1)=N_{\mu_1}(1)$.
Let us observe that this expression corresponds to the formula for the coefficients of 
$\eta_{\mu_1 \circlearrowright \mu_2}(n)$ given in [14] and 
quoted below for convenience, except that 
coefficients $N_{\mu_2}(k_1), \ldots , N_{\mu_2}(k_{r-1})$ are 
replaced by $N_{\nu_2}(k_1), \ldots , N_{\nu_2}(k_{r-1})$, respectively. 
Namely,
$$
N_{\mu_1 \circlearrowright\mu_2}(n)=
\sum_{r=2}^{n}N_{\mu_1}(r)\sum _{k_{1}+k_{2}+\ldots +k_{r}=n}
N_{\mu_2}(k_{1})\ldots
N_{\mu_2}(k_{r-1})N_{\mu_2}(k_r)
$$
for $n\geq 2$, whereas $N_{\mu_1 \circlearrowright\mu_2}(1)=N_{\mu_1}(1)$.
Let us select the indices in such a way that 
\begin{enumerate}
\item
$N_{\mu_{2}}(k_r)$ corresponds to the f-walks on the copy of ${\mathcal G}_{2}$ attached 
at the root $e_2$ to the root $e_1$ of ${\mathcal G}_{1}$,
\item 
the remaining coefficients $N_{\mu_{2}}(k_j)$, where $j=1, \ldots , r-1$, correspond to 
f-walks on the copies of ${\mathcal G}_{2}$ attached at the root $e_2$ to the 
remaining vertices of ${\mathcal G}_{1}$.
\end{enumerate}
Therefore, if we keep the first item above unchanged and we 
replace $\mu_2$ by $\nu_2$ in the second item above, this simply 
means that we count f-walks on the copy of ${\mathcal G}_{2}$ attached at the root $f_2$ to the $j$-th vertex of ${\mathcal G}_{1}$
for any $j=1, \ldots , r-1$. In other words, the walk counting 
is similar to that in the case of the comb product of graphs, except that 
all copies of ${\mathcal G}_{2}$ which are attached to any $v\in V_1^{\circ}$
are attached at the root $f_2$, not $e_2$. In order to use
the counting of alternating d-walks from $e$ to $e$ (as in the comb loop product, the alternating colors come from the loops of color different than the color of the preceding traversed edge), it suffices to add loops as follows:
\begin{enumerate}
\item loops of color $1$ to all vertices but the root $e_2$ of the copy of ${\mathcal G}_{2}$ attached to $e_1$ 
(as in the case of the comb product),
\item loops of color $1$ to all vertices but the root $f_2$ of each copy of ${\mathcal G}_{2}$ attached to the 
remaining vertices of ${\mathcal G}_{1}$.
\end{enumerate}
Observe that this is exactly the description of the c-comb loop product of graphs. Therefore, we conclude that 
the $\eta$-series of $\mu_1\circlearrowright_{\nu_2} \mu_2$ coincides with the series $\eta_{Z}$, where $Z=R_2R_1$, with $R_1$ and
$R_2$ are the adjacency matrices from the decomposition of Lemma 5.1 
of the c-comb loop product of graphs.
\hfill $\blacksquare$

\subsection*{\bf Acknowledgements}
This research has been partially supported by Narodowe Centrum Nauki (grant No. 2014/15/B/ST1/00166) and by the Wroc{\l}aw University of Science and Technology (project No. 0401/0121/17).

\end{document}